\documentclass[a4paper,11pt]{article}
\usepackage[nointlimits]{amsmath}
\usepackage{amssymb,amsmath,amsthm}
\usepackage{amsfonts}
\usepackage{latexsym}
\usepackage{graphicx}
\usepackage{color}
\usepackage{layout}
\usepackage[english]{babel}
\numberwithin{equation}{section} % to number equations within sections

\def\ds{\displaystyle}

\def\ZZ{\mathbb Z}

\def\bs{{\bf s}}
\def\bt{{\bf t}}
\def\br{{\bf r}}
\def\cP{{\cal P}}

\newcommand\stf{
\setlength{\unitlength}{.4pt}
\begin{picture}(40,20)
\put(10,2){\line(1,0){20}}
\put(10,2){\line(0,1){10}}
\put(20,2){\line(0,1){10}}
\put(15,7.5){\line(1,0){10}}
\put(30,2){\line(0,1){10}}
\end{picture}
}

\newtheorem{thm}{Theorem}[section]

\newtheorem{cor}[thm]{Corollary}

\begin{document}

\title{\bf New harmonic number identities with applications}
\author{{\sc Roberto Tauraso}\\
%\footnote{
Dipartimento di Matematica\\
Universit\`a di Roma ``Tor Vergata'', Italy\\
{\tt tauraso@mat.uniroma2.it}\\
{\tt http://www.mat.uniroma2.it/$\sim$tauraso}
%}
}

\date{}
\maketitle
%\thanks{}
\begin{abstract}
\noindent 
We determine the explicit formulas for the sum of products of homogeneous multiple
harmonic sums $\sum_{k=1}^n \prod_{j=1}^r H_k(\{1\}^{\lambda_j})$ when 
$\sum_{j=1}^r \lambda_j\leq 5$. 
We apply these identities to the study of two congruences modulo a power of a prime.
\end{abstract}

\makeatletter{\renewcommand*{\@makefnmark}{}
\makeatother}

%%%%%%%%%%%%%%%%%%%%%%%%%%%%%%%%%%%%%%%%%%%%%%%%%%%%%%%%%%%%%%%%%%%%%%%%%%
\section{Introduction} 
%%%%%%%%%%%%%%%%%%%%%%%%%%%%%%%%%%%%%%%%%%%%%%%%%%%%%%%%%%%%%%%%%%%%%%%%%%
Let $\bs=(s_1,s_2,\dots,s_d)$ be a vector whose entries are positive integers then we
define the {\sl multiple harmonic sum} (MHS for short) for $n\geq 0$ as
$$H_n(\bs)=\sum_{1\leq k_1<k_2<\cdots<k_d\leq n} \!{1\over k_1^{s_1}k_2^{s_2}\cdots k_d^{s_d}}.$$
We call $d$ and $|\bs|=\sum_{i=1}^d s_i$ its depth and its weight respectively.
This kind of sums has a long history, see for example \cite{Ho:07}, \cite{Zh:09} and their references.
 
In this note we present an algorithmic procedure to determine a closed formula for 
$$\sum_{k=1}^n \prod_{j=1}^m H_k(\bs_j)$$
involving products of MHS evaluated at $n$ and of total weight less or equal to 
$\sum_{j=1}^n|\bs_j|$.

\begin{itemize}

\item[({\bf i})]
The first step is to expand recursively any product $H_n(\bs)\cdot H_n(\bt)$ 
as a non-negative integral combination of MHS of weight $|\bs|+|\bt|$:
\begin{eqnarray*}H_n(\bs)\cdot H_n(\bt)&=&\sum_{\br\in \bs\stf \bt}H_n(\br)\\
&=&\sum_{\br\in \bs'\stf \bt}H_n(s_1,\br)+\sum_{\br\in \bs\stf \bt'}H_n(v_1,\br)+
\sum_{\br\in \bs'\stf \bt'}H_n(s_1+t_1,\br)
\end{eqnarray*}
where $\bs\stf \bt$ is the so-called {\sl stuffle product} (or {\sl quasi-shuffle product}) 
of $\bs=(s_1, \bs')$ and $\bt=(t_1,\bt')$ (note that if $\bt$ is empty then $\bs\stf \bt=\bs$ and
$\sum_{\br\in \bs\stf \bt}H_n(\br)$ is simply $H_n(\bs)$).

\item[({\bf ii})] 
The second step is to sum over $k$ any MHS by using recursively the following rule:
\begin{eqnarray*}
\sum_{k=1}^n H_k(\bs)&=&\sum_{k=1}^n\sum_{j=1}^k {H_{j-1}(\bs')\over j^{s_d}}
=\sum_{j=1}^n {H_{j-1}(\bs')\over j^{s_d}}\sum_{k=j}^n 1
=\sum_{j=1}^n {H_{j-1}(\bs')\over j^{s_d}}(n+1-j)\\
&=&
(n+1)H(\bs)-
\left\{\begin{array}{ll}
H_n(\bs',s_d-1) &\mbox{if $s_d>1$}\\
\ds \sum_{k=1}^{n} H_k(\bs') -H_n(\bs')&\mbox{if $s_d=1$}
\end{array}
\right..
\end{eqnarray*}
where $\bs=(\bs',s_d)$.
It is easy to see that the MHS involved in the final formula have weight less or equal to $|\bs|$.

\item[({\bf iii})] The last step is to simplify the formula by collecting terms by using the stuffle product relations. This can be done in various ways and the choices depend on the application of the identity.
\end{itemize}

This procedure will be employed in the next section 
and it will generate a large bunch of identities.
We give two applications of these identities (we hope that many more will come out in the future). The first is about the sum 
$$\sum_{k=0}^{n}(-1)^{ak}{n\choose k}^a$$
which seems to have a closed form only for a very few values of the parameter $a\in\ZZ$, namely for $-1\leq a\leq 3$.
In \cite{CG:02}, Cai and Granville studied a related congruence by showing that for any prime $p\geq 5$
$$\sum_{k=0}^{p-1}(-1)^{ak}{p-1\choose k}^a\equiv {ap-2\choose p-1}\pmod{p^4}.$$
Here we present the following extension (see also \cite{Pa:09} and \cite{CD:06} for similar results): 
 
\begin{thm} Let $p>5$ be a prime then for any $a\in\ZZ$
$$\sum_{k=0}^{p-1}(-1)^{ak}{p-1\choose k}^a\equiv
{(a-1)p\over ap-1}\left(1+{a(a+1)(3a-2)\over 6}p^3 X_p\right) \pmod{p^6}$$
where $X_p={B_{p-3}\over p-3}-{B_{2p-4}\over 4p-8}$ and $B_n$ denotes the $n$-th Bernoulli number.
\end{thm}

Moreover, by using the previous theorem for $a=-2$ we get 

\begin{cor} Let $p>5$ be a prime then
$$\sum_{k=1}^{p-1}{1\over k}{2k\choose k}\equiv -{16\over 3}p^2 X_p \pmod{p^4}.$$
\end{cor}

which improves the (mod $p^3$) result contained in \cite{ST:08}.

%%%%%%%%%%%%%%%%%%%%%%%%%%%%%%%%%%%%%%%%%%%%%%%%%%%%%%%%%%%%%%%%%%%%%%%%%%
\section{MHS: identities} 
%%%%%%%%%%%%%%%%%%%%%%%%%%%%%%%%%%%%%%%%%%%%%%%%%%%%%%%%%%%%%%%%%%%%%%%%%%

By following the procedure introduced in the previous section, 
we found an explicit formula for any
sum of products of homogeneous MHS like $H_k(\{1\}^d)$ up to order $5$
($\{1\}^d$ means that the number $1$ is repeated $d$ times):
\begin{eqnarray*}
&&\sum_{k=1}^n H_k(1)=(n+1)H_n(1)-n\,,\\
&&\sum_{k=1}^n H_k(\{1\}^2)=(n+1)H_n(\{1\}^2)-nH_n(1)+n\,,\\
&&\sum_{k=1}^n H_k^2(1)=(n+1)H_n^2(1)-(2n+1)H_n(1)+2n\,,\\
&&\sum_{k=1}^n H_k(\{1\}^3)=(n+1)H_n(\{1\}^3)+n\left(H_n(1)-{1\over 2}H^2_n(1)\right)+{n\over 2}H_n(2)-n\,,\\
&&\sum_{k=1}^n H_k(1)H_k(\{1\}^2)=(n+1)H_n(1)H_n(\{1\}^2)+(3n+1)\left(H_n(1)-{1\over 2}H^2_n(1)\right)+{n+1\over 2}H_n(2)-3n\,,\\
&&\sum_{k=1}^n H_k^3(1)=(n+1)H_n^3(1)+(6n+3)\left(H_n(1)-{1\over 2}H^2_n(1)\right)+{1\over 2}H_n(2)-6n.
\end{eqnarray*}
It's interesting to note that the formulas for $\sum_{k=1}^n H_k^r(1)$ when $r=1,2,3$ 
appear as Entry 8 at page 94 in \cite{Be:98}.
To illustrate the procedure we show how to obtain $\sum_{k=1}^n H_k(1)H_k(\{1\}^2)$.

\noindent By (i)
$$H_k(1)H_k(\{1\}^2)=3H_k(\{1\}^3)+H_k(2,1)+H_k(1,2).$$
By (ii)
\begin{eqnarray*}
&&\sum_{k=1}^n H_k(\{1\}^3)=(n+1)H_n(\{1\}^3)-nH_n(\{1\}^2)+nH_n(1)-n\,,\\
&&\sum_{k=1}^n H_k(2,1)=(n+1)H_n(2,1)-nH_n(2)+H_n(1)\,,\\
&&\sum_{k=1}^n H_k(1,2)=(n+1)H_n(1,2)-H_n(\{1\}^2).
\end{eqnarray*}
Finally by (iii), since by (i)
$$H_n(1)H_n(\{1\}^2)=3H_n(\{1\}^3)+H_n(2,1)+H_n(1,2)\quad\mbox{and}\quad 2H_n(\{1\}^2)=H_n^2(1)-H_n(2),$$
we get the formula given above.

\noindent The formulas when the total weight is $4$ are contained in the next table: 
the sum $$\sum_{k=1}^n f_k-(n+1)f_n,$$
 where $f_n$ is an entry in the first row,
is equal to the sum of the entries of the first column each multiplied by the linear polynomial $an+b$ contained in the intersection
of the chosen row and column.

\begin{center}
\begin{tabular}{|c|c|c|c|c|c|}\hline
$\sum_{k=1}^n f_k-(n+1)f_n$&$\ds H_n(\{1\}^4)$ &$\ds H_n^2(\{1\}^2)$ &$\ds H_n(1)H_n(\{1\}^3)$ &$\ds H^2_n(1)H_n(\{1\}^2)$ &$\ds H^4_n(1)$\\\hline
$\sum_{k=1}^3{(-1)^{k-1}\over k!}H_n^k(1)$ &
$ -n $ &$ -6n-2 $ &$ -4n-1$ &$ -12n-5$ &$ -24n-12$\\\hline
$ {1\over 2}H_n(2)$ &$ -n $ &$-2n-2 $ &$-2n-1$ &$-2n-3$ &$-4$\\\hline
$ {1\over 3}H_n(3)$ &$ -n $ &$1 $ &$-n-1$ &$1$ &$3$\\\hline
$ {1\over 2}H_n(1)H_n(2)$ &$ n $ &$2n$ &$2n+1$ &$2n+1$ &$0$\\\hline
$ H_n(1,2)$ &$ 0$ &$1$ &$0$ &$1$ &$2$\\\hline
$ n$ &$ 1$ &$6$ &$4$ &$12$ &$24$\\\hline
\end{tabular}
\end{center}

\noindent The large table which gives the formulas when the total weight is $5$ is in the Appendix.

%%%%%%%%%%%%%%%%%%%%%%%%%%%%%%%%%%%%%%%%%%%%%%%%%%%%%%%%%%%%%%%%%%%%%%%%%%
\section{MHS: congruences} 
%%%%%%%%%%%%%%%%%%%%%%%%%%%%%%%%%%%%%%%%%%%%%%%%%%%%%%%%%%%%%%%%%%%%%%%%%%

Among the various known results about MHS modulo power of a prime,
the following ones will be crucial for us: for any prime $p>5$
$$\begin{array}{lllll}
H_{p-1}(1)\equiv 2p^2X_p &\pmod{p^4}\,,\\
H_{p-1}(2)\equiv -4pX_p &\pmod{p^3}\,,\\
H_{p-1}(3)\equiv 0 &\pmod{p^2}\,,\\
H_{p-1}(1,2)\equiv -6X_p &\pmod{p^2}\,,\\
H_{p-1}(4) \equiv H_{p-1}(\{1\}^2,2)\equiv H_{p-1}(1,3) \equiv 0  &\pmod{p}.
\end{array}
$$
where $X_p={B_{p-3}\over p-3}-{B_{2p-4}\over 4p-8}$ and $B_n$ denotes the $n$-th Bernoulli number (see \cite{Suzh:00} for the MHS of depth $1$,
see \cite{Zh:09} and \cite{Ho:07} for all the MHS of depth $>1$ with the exception of   $H_{p-1}(1,2)$ modulo $p^2$ which has been established in \cite{Ta:09}).

\noindent Note that every homogeneous MHS can be expressed in terms of MHS of depth $1$.
More precisely (see Theorem 2.3 in \cite{Ho:07}): given a positive integer $d$ then
for any unordered partition $\lambda=\{\lambda_1,\lambda_2,\dots,\lambda_r\}$ of $d$ there is
some integer $c_{\lambda}$ such that 
$$d!H_n(\{1\}^d)=\sum_{\lambda\in P(d)}c_{\lambda}\prod_{i=1}^r H_n(\lambda_i).$$
For example:
\begin{eqnarray*}
&&2H_n(\{1\}^2)=H_n^2(1)-H_n(2)\,,\\
&&6H_n(\{1\}^3)=H_n^3(1)-3H_n(1)H_n(2)+2H_n(3)\,,\\
&&24H_n(\{1\}^4)=H_n^4(1)-6H_n^2(1)H_n(2)+8H_n(1)H_n(3)+3H_n^2(2)-6H_n(4).
\end{eqnarray*}
Hence, for any prime $p>5$
$$\begin{array}{lllll}
H_{p-1}(\{1\}^2)\equiv  2pX_p &\pmod{p^3}\,,\\
H_{p-1}(\{1\}^3) \equiv 0 &\pmod{p^2}\,,\\
H_{p-1}(\{1\}^4) \equiv 0  &\pmod{p}.
\end{array}
$$

\noindent Therefore, by the previous identities and congruences, it follows that for any prime $p>5$:

\begin{center}
\begin{tabular}{|ll|}\hline
$\sum_{k=1}^{p-1} H_k(1)\equiv -(p-1)+2p^3X_p$ &$\pmod{p^5}$\\\hline
$\sum_{k=1}^{p-1} H_k(\{1\}^2)\equiv (p-1)+(4-2p)p^2X_p$ &$\pmod{p^4}$\\
$\sum_{k=1}^{p-1} H_k^2(1)\equiv 2(p-1)+(2-4p)p^2X_p$ &$\pmod{p^4}$\\\hline
$\sum_{k=1}^{p-1}H_k(\{1\}^3)\equiv -(p-1)+(2-4p)pX_p$ &$\pmod{p^3}$\\
$\sum_{k=1}^{p-1}H_k(1)H_k(\{1\}^2)\equiv -3(p-1)+(-6p)pX_p$ &$\pmod{p^3}$\\
$\sum_{k=1}^{p-1}H_k^3(1)\equiv -6(p-1)+(-2-6p)pX_p$ &$\pmod{p^3}$\\\hline
$\sum_{k=1}^{p-1}H_k(\{1\}^4)\equiv (p-1)+(-2p)X_p$ &$\pmod{p^2}$\\
$\sum_{k=1}^{p-1}H_k^2(\{1\}^2)\equiv 6(p-1)+(-6)X_p$ &$\pmod{p^2}$\\
$\sum_{k=1}^{p-1}H_k(1)H_k(\{1\}^3)\equiv 4(p-1)+(-2p)X_p$ &$\pmod{p^2}$\\
$\sum_{k=1}^{p-1}H_k^2(1)H_k(\{1\}^2)\equiv 12(p-1)+(-6+2p)X_p$ &$\pmod{p^2}$\\
$\sum_{k=1}^{p-1}H_k^4(1)\equiv 24(p-1)+(-12+8p)X_p$ &$\pmod{p^2}$\\\hline
$\sum_{k=1}^{p-1}H_k(\{1\}^5)\equiv 1$ &$\pmod{p}$\\
$\sum_{k=1}^{p-1}H_k(1)H_k(\{1\}^4)\equiv 5$ &$\pmod{p}$\\
$\sum_{k=1}^{p-1}H_k(\{1\}^2)H_k(\{1\}^3)\equiv 10+6X_p$ &$\pmod{p}$\\
$\sum_{k=1}^{p-1}H_k^2(1)H_k(\{1\}^3)\equiv 20+6X_p$ &$\pmod{p}$\\
$\sum_{k=1}^{p-1}H_k(1)H_k^2(\{1\}^2)\equiv 30+18X_p$ &$\pmod{p}$\\
$\sum_{k=1}^{p-1}H_k^3(1)H_k(\{1\}^2)\equiv 60+30X_p$ &$\pmod{p}$\\
$\sum_{k=1}^{p-1}H_k^5(1)\equiv 120+60X_p$ &$\pmod{p}$\\\hline
\end{tabular}
\end{center}
Note that $\sum_{k=1}^{p-1} H_k^r(1)$ (mod $p^{4-r}$) for $r=1,2,3$ have been established by Z. W. Sun in \cite{Suzw:09}.

\section{Proof of Theorem 1.1 and Corollary 1.2} 

\begin{proof}[Proof of Theorem 1.1] 
Assume that $p>5$ is a prime, then for $k=1,\dots, p-1$ we have that
$$(-1)^k{p-1\choose k}=\prod_{j=1}^{k}\left(1-{p\over j}\right)\equiv1+\sum_{j=1}^5 (-p)^j H_k(\{1\}^j)\pmod{p^6}.$$
Hence
$$(-1)^{ak}{p-1\choose k}^a\equiv 1+\sum_{j=1}^5(-p)^j
\sum_{r=1}^j {a \choose r} \sum_{\lambda\in \cP(j,r)} {j \choose \lambda_1,\dots,\lambda_r}
\prod_{i=1}^r H_k(\{1\}^{\lambda_i})\pmod{p^6}$$
where $\cP(j,r)$ is the set of the integer partitions $\lambda$ of $j$ into $r$ parts.

\noindent By summing over $k$ we find
$$\sum_{k=0}^{p-1}(-1)^{ak}{p-1\choose k}^a=p+\sum_{j=1}^5(-p)^j
\sum_{r=1}^j {a \choose r} \sum_{\lambda\in \cP(j,r)} {j \choose \lambda_1,\dots,\lambda_r}
\sum_{k=1}^{p-1}\prod_{i=1}^r H_k(\{1\}^{\lambda_i})\pmod{p^6}.$$
Finally, by the congruences established in the previous section we can compute
$$(-p)^j\sum_{k=1}^{p-1}\prod_{i=1}^r H_k(\{1\}^{\lambda_i})\pmod{p^{6}}$$
for any partition $\lambda \in\cP(j,r)$ and we get easily the result.
\end{proof} 

\begin{proof}[Proof of Corollary 1.2] 
In \cite{St:47} Staver proved that for any integer $n\geq 1$
$$\sum_{k=1}^{n}{1\over k}{2k\choose k}={2n\choose n}{2n+1\over 3n^2}\sum_{k=0}^{n-1}{n-1\choose k}^{-2}.$$
Letting $n=p$, by Theorem 1.1 for $a=-2$ we have
\begin{eqnarray*}
\sum_{k=1}^{p-1}{1\over k}{2k\choose k}&=&
{1\over p}{2p\choose p}\left({2p+1\over 3p}\sum_{k=0}^{p-1}{p-1\choose k}^{-2}-1\right)\\
&\equiv& {2\over p}{2p-1\choose p-1}\left(\left(1-{8\over 3}p^3 X_p\right)-1\right)\equiv -{16\over 3}p^2 X_p\pmod{p^4}
\end{eqnarray*}
where in the last step we used the fact that ${2p-1\choose p-1}\equiv 1 \pmod{p^3}$ by Wolstenholme theorem.
\end{proof}

%%%%%%%%%%%%%%%%%%%%%%%%%%%%%%%%%%%%%%%%%%%%%%%%%%%%%%%%%%%%%%%%%%%%%%%

%%%%%%%%%%%%%%%%%%%%%%%%%%%%%%%%%%%%%%%%%%%%%%%%%%%%%%%%%%%%%%%%%%%%%%%%%%
\section*{Appendix} 
%%%%%%%%%%%%%%%%%%%%%%%%%%%%%%%%%%%%%%%%%%%%%%%%%%%%%%%%%%%%%%%%%%%%%%%%%%

\begin{center}
\begin{tabular}{|c|c|c|c|c|}\hline
$\sum_{k=1}^n f_k-(n+1)f_n$&$\ds H_n(\{1\}^5)$ &$\ds H_n(1)H_n(\{1\}^4)$ &$\ds H_n(1,1)H_n(\{1\}^3)$ &$\ds H^2_n(1)H_n(\{1\}^3)$ \\\hline
$\sum_{k=1}^4{(-1)^{k-1}\over k!}H_n^k(1)$ &
$ n $ &$ 5n+1 $ &$ 10n+3$ &$ 20n+7$\\\hline
$ {1\over 2}H_n(2)$ &$ n $ &$ 3n+1 $ &$4n+3$ &$6n+5$\\\hline
$ {1\over 3}H_n(3)$ &$ n $ &$2n+1 $ &$n$ &$2n+1$\\\hline
${1\over 2} H_n(1)H_n(2)$ &$ -n $ &$-3n-1$ &$-4n-1$ &$-6n-3$\\\hline
$ H_n(1,2)$ &$ 0$ &$0$ &$-1$ &$-1$\\\hline
$ {1\over 4}H_n(4)$ &$ n $ &$n+1 $ &$-1$ &$-1$\\\hline
$ {1\over 8}H^2_n(2)$ &$ -n $ &$-(n+1) $ &$-(2n-1)$ &$1$\\\hline
$ {1\over 4}H_n^2(1)H_n(2)$ &$ n $ &$3n+1 $ &$4n+1$ &$6n+3$\\\hline
$ {1\over 3}H_n(1)H_n(3)$ &$ -n $ &$-2n-1 $ &$-n$ &$-2n-1$\\\hline
$ H_n(1,3)$ &$ 0$ &$0$ &$0$ &$0$\\\hline
$ H_n(1,1,2)$ &$ 0$ &$0$ &$1$ &$1$\\\hline
$ n$ &$ -1$ &$-5$ &$-10$ &$-20$\\\hline
\end{tabular}

\vspace*{5mm}

\begin{tabular}{|c|c|c|c|c|}\hline
$\sum_{k=1}^n f_k-(n+1)f_n$&$\ds H_n(1)H_n^2(\{1\}^2)$ &$\ds H_n^3(1)H_n(\{1\}^2)$ &$\ds H_n^5(1)$\\\hline
$\sum_{k=1}^4{(-1)^{k-1}\over k!}H_n^k(1)$ &
$ 30n+12 $ &$ 60n+27 $ &$ 120n+60$\\\hline
$ {1\over 2}H_n(2)$ &$ 6n+8 $ &$ 6n+13 $ &$20$\\\hline
$ {1\over 3}H_n(3)$ &$ -3 $ &$-6 $ &$-15$\\\hline
${1\over 2} H_n(1)H_n(2)$ &$ -6n-2 $ &$-6n-3$ &$0$\\\hline
$ H_n(1,2)$ &$ -3$ &$-5$ &$-10$\\\hline
$ {1\over 4}H_n(4)$ &$-2 $ &$-3$ &$-4$\\\hline
$ {1\over 8}H^2_n(2)$ &$ -2n+4 $ &$9 $ &$20$\\\hline
$ {1\over 4}H_n^2(1)H_n(2)$ &$6n+2 $ &$6n+3 $ &$0$\\\hline
$ {1\over 3}H_n(1)H_n(3)$ &$ 0 $ &$0 $ &$0$\\\hline
$ H_n(1,3)$ &$1$ &$2$ &$5$\\\hline
$ H_n(\{1\}^2,2)$ &$3$ &$5$ &$10$\\\hline
$ n$ &$-30$ &$-60$ &$-120$\\\hline
\end{tabular}
\end{center}
\end{document}